\documentclass{amsart}
\usepackage{amsmath,amsfonts,amssymb,amsthm,epsfig}

\newtheorem{Theorem}{Theorem}[section]
\newtheorem{Proposition}[Theorem]{Proposition} 
\newtheorem{Lemma}[Theorem]{Lemma}

\newtheorem{Question}{Question}

\newcommand{\wt}{\operatorname{wt}}
\newcommand{\wi}{\mathbf{i}}
\newcommand{\bigdot}{\bullet}
\newcommand{\rev}{\mathrm{rev}}

\theoremstyle{definition}

\newtheorem{Remark}[Theorem]{Remark}

\begin{document}

\title[A definition of the crystal commutor using Kashiwara's involution]{A definition of the crystal commutor using Kashiwara's involution}

\author{Joel Kamnitzer}
\email{jkamnitz@aimath.org}
\address{American Institute of Mathematics\\ Palo Alto, CA}

\author{Peter Tingley}
\email{pwtingle@math.berkeley.edu}
\address{UC Berkeley, Department of Mathematics\\ Berkeley, CA}

 \thanks{The second author was supported by the RTG grant DMS-0354321}

\begin{abstract}
Henriques and Kamnitzer defined and studied a commutor for the category of crystals of a finite dimensional complex reductive Lie algebra. We show that the action of this commutor on highest weight elements can be expressed very simply using Kashiwara's involution on the Verma crystal.

\end{abstract}

\date{\today}
\maketitle
 
\section{Introduction}
Let $ \mathfrak{g} $ be a complex reductive Lie algebra. If $ A $ and $ B $ are crystals of representations of $\mathfrak{g}$, then $ A \otimes B $ and $ B \otimes A $ are isomorphic.  However the map $ (a,b) \mapsto (b,a) $ is not an isomorphism.  In \cite{HK}, following an idea of Berenstein, A. Henriques and the first author construct an explicit isomorphism $ \sigma_{A,B}: A \otimes B \rightarrow B \otimes A $, which they call the commutor.  This commutor is involutive and satisfies the cactus relation, a certain axiom involving triple tensor products (see section \ref{questions}).

Consider the following alternative definition for a commutor.  First notice that we only need to define $\sigma_{A,B}$ when $A$ and $B$ are irreducible. Also, a crystal isomorphism is uniquely defined by the images of highest weight elements. So, for each highest weight element $b_\lambda \otimes c \in B_\lambda \otimes B_\mu$, we need to specify its image $b_\mu \otimes b \in B_\mu \otimes B_\lambda$. We do this using Kashiwara's involution $*$ on $B_\infty$. By the properties of $*$, if 
 $b_\lambda \otimes c$ is a highest weight element in $B_\lambda \otimes B_\mu$, then $*c \in B_\lambda$ (where we identify $B_\lambda$ and $B_\mu$ with their images in $B_\infty$), and $b_\mu \otimes *c$ is highest weight. Therefore we can define a crystal commutor by specifying that each highest weight element $b_\lambda \otimes c$ is taken to $b_\mu \otimes *c$. In this note we show that this definition gives the same commutor as that studied by Henriques and Kamnitzer. 

The original definition of the commutor used the Sch\"utzenberger involution on each $ B_\lambda $, while this definition uses Kashiwara's involution on $ B_\infty $.  Thus one way to interpret our result is that it gives a non-trivial relationship between these two involutions.  The Sch\"utzenberger involution does not exist for crystals of non-finite symmetrizable Kac-Moody Lie algebras, however Kashiwara's involution does.  Hence this work extends the definition of the commutor to highest weight crystals of symmetrizable Kac-Moody Lie algebras.

\subsection{Acknowledgments}
We thank A. Berenstein, M. Haiman, and M. Vazriani for useful discussions and the referee for helpful comments.  The first author thanks the American Institute of Mathematics for support and UC Berkeley for hospitality.

\section{Background}

\subsection{Notation and terminology} \label{notation} 
We include only a brief review of some basic facts about crystals. For the most part we follow the conventions from the review article \cite{Kas2}, which we recommend for a more detailed overview of the subject.
\begin{itemize}
\item Let $ \mathfrak{g} $ be a complex reductive Lie algebra. 

\item Let $I $ denote the set of vertices of the Dynkin diagram of $\mathfrak{g}$.  

\item Let $ \{ \alpha_i \}_{i \in I}, \{\alpha^\vee_i \}_{i \in I} $ denote the positive roots and coroots of $\mathfrak{g}$. 

\item Let $ \{s_i\}_{i \in I} $ denote the generators of the Weyl group.

\item Let $ w_0 $ denote the long element of the Weyl group. 

\item Let $\langle \cdot, \cdot \rangle$ denote the pairing between weight space and coweight space.

\item  Let $\Lambda$ denote the set of weights of $\mathfrak{g}$, $ \{\Lambda_i\}_{i \in I} $ the set of fundamental weights, and  $ \Lambda_+ $ the set of dominant weights.  

\item A crystal for $\mathfrak{g} $ is a finite set $ B $ along with maps $ e_i, f_i : B \rightarrow B \sqcup 0$ for each $ i \in I $, and a map $ \mathrm{wt} : B \rightarrow \Lambda $, satisfying a certain set of axioms.  These axioms may be found in \cite[Section 7.2]{Kas2}.

\item For $ \lambda \in \Lambda_+ $, let $B_\lambda$ denote the crystal corresponding to the irreducible representation $V_\lambda$ of $\mathfrak{g}$.  

\item An element of a crystal is called highest (resp. lowest) weight if it is killed by all $ e_i$ (resp. all $f_i$).  We use  $b_\lambda$ and $b_\lambda^{low}$ to denote the unique highest and lowest  weight elements of $ B_\lambda $.

\item For any crystal $B$ and any $b \in B$, let  $\varepsilon_i(b) = \max \{ n : e_i^n (b) \neq 0 \}$. Let $\varepsilon(b) \in \Lambda_+$ be the unique weight such that, for all $i \in I$, $(\varepsilon(b), \alpha_i^\vee)= \varepsilon_i(b)$.

\item Similarly, let  $\varphi_i (b) = \max \{ n : f_i^n (b) \neq 0 \}$ and $\varphi(b) \in \Lambda_+$ be the unique weight such that, for all $i \in I$, $(\varphi(b), \alpha_i^\vee= \varphi_i(b)$.

\item The weight of $b \in B$ is $\wt(b) = \varphi(b)-\varepsilon(b)$.

\item There is a tensor product rule for crystals corresponding to the tensor product for representations of $\mathfrak{g}$. The underlying set of $ A \otimes B $ is $ A \times B$ (whose elements we denote $ a \otimes b$) and the actions of $ e_i $ and $ f_i $ are given by the following rules:

\begin{equation} \label{edef} e_i (a \otimes b)=
\begin{cases}
e_i  (a) \otimes b, \quad \text{if}\quad  \varphi_i(a) \geq \varepsilon_i(b)\\
a \otimes e_i  (b),\quad \text{otherwise}
\end{cases} \end{equation}
\begin{equation} \label{fdef} f_i (a \otimes b)=
\begin{cases}
f_i  (a) \otimes b, \quad \text{if} \quad  \varphi_i(a) > \varepsilon_i(b)\\
a \otimes f_i (b),\quad \text{otherwise}.
\end{cases} \end{equation}

\end{itemize}
\subsection{Kashiwara's involution on $B_\infty$} \label{kinv}
For any dominant weights $\gamma$ and $\lambda$, there is an inclusion of crystals $ B_{\gamma+\lambda} \rightarrow B_\lambda \otimes B_\gamma$ which sends $b_{\gamma+\lambda}$ to $b_\gamma \otimes b_\lambda$. The following is immediate from the tensor product rule:

\begin{Lemma} \label{in_gamma}
The image of the inclusion $ B_{\lambda+\gamma} \rightarrow B_\gamma \otimes B_\lambda $ contains all elements of the form $ b \otimes b_\gamma $ for $ b \in B_\lambda $.  \qed
\end{Lemma}

Lemma \ref{in_gamma} defines a map $ \iota_{\lambda}^{\lambda + \gamma} : B_\lambda \rightarrow B_{\lambda + \gamma} $ which is $ e_i $ equivariant and takes $ b_\lambda $ to $ b_{\lambda + \gamma} $.  
These maps make $\{ B_\lambda \}$ into a directed system, and the limit of this system is $ B_\infty $.  There are $ e_i $ equivariant maps $ \iota_\lambda^\infty : B_\lambda \rightarrow B_\infty $. When there is no danger of confusion we denote $\iota_\lambda^\infty$ simply by $\iota$.  

The infinite set $ B_\infty $ has additional combinatorial structure. In particular, we will need:
\begin{enumerate}
\item The map $\tau : B_\infty \rightarrow \Lambda_+ $ defined by  $\tau(b) = \min \{ \lambda :  b \in \iota (B_\lambda) \}.$
\item The map  $ \varepsilon : B_\infty \rightarrow \Lambda_+ $ given by, for any $b \in B_\infty$ and any $\lambda$ such that $b \in \iota(B_\lambda)$, $\varepsilon(b)= \varepsilon(\iota^{-1}(b))$, where $\varepsilon$ is defined on $B_\lambda$ as in Section \ref{notation}. 
\item Kashiwara's involution $*$ (for the construction of this involution see \cite[Theorem 2.1.1]{Kas1}).
\end{enumerate}

These maps are related by the following result of Kashiwara.

\begin{Proposition}[\cite{Kas2}, Prop. 8.2] \label{th:prop}
Kashiwara's involution preserves weights and satisfies 
\begin{equation*}
\tau(*b) = \varepsilon(b), \quad \varepsilon(*b) = \tau(b). 
\end{equation*} 
\end{Proposition}

\begin{Remark}
All of this combinatorial structure can be seen easily using the MV polytope model \cite{Kam}.  The inclusions $ \iota $ correspond to translating polytopes.  The maps $ \tau $ and $\varepsilon $ are given by counting the lengths of edges coming out of the top and bottom vertices.  The involution $ * $ corresponds to negating a polytope.  From this description, the proof of the above proposition is immediate.
\end{Remark}

\subsection{The commutor} \label{commutor_section}

We now recall the definition of the commutor from \cite[Section 2.2]{HK}.  Let $ \theta : I \rightarrow I $ be the involution such that $ -w_0 \cdot \alpha_i = \alpha_{\theta(i)} $.  Recall that each crystal $ B_\lambda $ comes with an involution $ \xi_\lambda $ which acts by $ w_0 $ on weights and exchanges the action of $ e_i $ and $f_{\theta(i)} $.  These involutions can be extended to a map $ \xi_B : B \rightarrow B $ for any crystal $ B $ and they lead to the definition of the commutor for crystals.  Namely,
\begin{align} \label{commuter_definition}
 \nonumber \sigma_{B,C} : B\otimes C &\rightarrow C \otimes B \\
 b \otimes c &\mapsto \xi_{C \otimes B}(\xi_C(c) \otimes \xi_B(b)) = \mathrm{Flip} \circ \xi_B \otimes \xi_C (\xi_{B \otimes C} (b \otimes c)).
 \end{align}
The second expression here is just the inverse of the first expression, and the equality is proved in \cite[Proposition 2]{HK}.

\section{Main theorem} \label{mainsection}

A crystal isomorphism $B_\lambda \otimes B_\mu \rightarrow B_\mu \otimes B_\lambda$ is uniquely defined by the images of the highest weight elements in
$B_\lambda \otimes B_\mu$. These are all of the form $b_\lambda \otimes c$, and must be sent to highest
weight elements of $B_\mu \otimes B_\lambda$, which in turn are of the form $b_\mu \otimes b$. It follows from the tensor product rule for crystal that $b_\lambda \otimes c$ is highest weight in $B_\lambda \otimes B_\mu $ if and only if $ \varepsilon(c) \le \lambda$.

As in Section \ref{kinv}, $B_\lambda$ and $B_\mu$ embed in $B_\infty$.   Let $ b_\lambda \otimes c $ be a highest weight element in $ B_\lambda \otimes B_\mu $.
 Since $\varepsilon (c) \leq \lambda$, by Proposition \ref{th:prop}, $\tau(*c) \leq \lambda$, or equivalently $*c \in \iota (B_\lambda)$. For this reason $*c$ can be considered an element of $B_\lambda$.  Also  $ \tau(c) \le \mu$, which implies that $\varepsilon(*c) \le \mu $, and so $ b_\mu \otimes *c $ is highest weight. So there is a unique isomorphism of crystals $B_\lambda \otimes B_\mu \rightarrow B_\mu \otimes B_\lambda$ which takes each highest weight element $b_\lambda \otimes c$ to $b_\mu \otimes *c$. The following shows that this isomorphism is equal to the crystal commutor.

\begin{Theorem}
\label{maintheorem}
If $ b_\lambda \otimes c$ is a highest weight element in $B_\lambda \otimes B_\mu$, then
$\sigma_{B_\lambda,B_\mu} (b_\lambda \otimes c)= b_\mu \otimes *c$.

\end{Theorem}

 \section{Proof}
 \label{proof}

One of the main tools we will need is the notion of Kashiwara data (also called string data), first studied by Kashiwara (see for example \cite{Kas2} section 8.2).
Fix a reduced word ${\bf i}$ for $w_0$, by which we mean $ \wi = (i_1, \dots, i_m) $, where each $ i_k $ is a node of the Dynkin diagram, and $ w_0 = s_{i_1} \cdots s_{i_m} $. The downward Kashiwara data for $b \in B_\lambda$ with respect to $ \wi $ is the sequence of non-negative integers $(p_1, \ldots p_m)$ defined by
\begin{gather*}
p_1 := \varphi_{i_1}(b), \quad p_2 := \varphi_{i_2}(f_{i_1}^{p_1} b), \quad \dots, \quad p_m := \varphi_{i_m}(f_{i_{m-1}}^{p-1} \ldots f_{i_1}^{p_1} b).
\end{gather*}
That is, we apply the lowering operators in the direction of $i_1$ as far as we can, then in the direction $i_2$, and so on. The following result is due to Littelmann \cite[section 1]{L}.

\begin{Lemma} \label{th:bottom} After we apply these steps, we reach the lowest element of the crystal. That is:
$$f_{i_m}^{p_m} \ldots f_{i_1}^{p_1} b = b_\lambda^{low}.$$

Moreover, the map $B_\lambda \rightarrow \mathbb{N}^m$ taking $b \rightarrow (p_1, \ldots p_m)$
is injective. \qed
\end{Lemma}

Similarly, the upwards Kashiwara data for $b \in B_\lambda$ with respect to ${\bf i}$ is the sequence $(q_1, \ldots q_m)$ defined by
\begin{gather*}
q_1 := \varepsilon_{i_1}(b), \quad q_2 := \varepsilon_{i_2}(e_{i_1}^{q_1} b), \quad \dots, \quad q_m := \varepsilon_{i_m}(e_{i_{m-1}}^{q-1} \ldots e_{i_1}^{q_1} b).
\end{gather*}

We introduce the notation $ w_k^\wi := s_{i_1} \cdots s_{i_k} $.

\begin{Lemma} \label{KDL}  In the crystal $ B_\lambda $, we have the following:
\begin{enumerate}
\item The downward Kashiwara data for $b_\lambda$ is given by $p_k =\langle w_{k-1}^\wi \cdot \alpha_{i_k}^\vee , \lambda \rangle.$ 
\item For each $k$, $\varepsilon_{i_k}  ( f_{i_{k-1}}^{p_{k-1}}  \ldots f_{i_1}^{p_1} b_\lambda )= 0.$
\end{enumerate}
\end{Lemma}

\begin{proof} Let $(p_1, \ldots p_m)$ be the downwards Kashiwara data for $b_\lambda$, and
let $\mu_k$ be the weight of $ f_{i_k}^{p_k} \ldots f_{i_1}^{p_1} b_\lambda$. Since $f_{i_k}^{p_k} \ldots f_{i_1}^{p_1} b_\lambda$ is the end of an $\alpha_{i_k}$ root string, we see that
\begin{equation} \label{reflect} \mu_k= s_{i_k} \cdot \mu_{k-1} -  a_k\alpha_{i_k}, \end{equation} 
 where $a_k = \varepsilon_{i_k}  ( f_{i_{k-1}}^{p_{k-1}}  \ldots f_{i_1}^{p_1} b_\lambda )$.  Using this fact at each step,
$$\mu_m = w_0 \cdot \lambda - \sum_{k=1}^m a_k s_{i_m} \ldots s_{i_{k+1}} \cdot \alpha_{i_k}.$$
By Lemma \ref{th:bottom}, we know that $f_{i_m}^{p_m} \ldots f_{i_1}^{p_1} b_\lambda = b_\lambda^{low}$, so that $\mu_m =  w_0 \cdot \lambda$. 
Hence 
$$ \sum_{k=1}^m a_k s_{i_m} \ldots s_{i_{k+1}} \cdot \alpha_{i_k} = 0. $$
Now, $s_{i_m} \cdots s_{i_k} $ is a reduced word for each $k$, which implies that
$s_{i_m} \ldots s_{i_{k+1}} \alpha_{i_k}$ is a positive root. Thus
each $a_k$ is zero, proving part (ii).

Equation (\ref{reflect}) now shows that $\mu_k = s_{i_k} \ldots s_{i_1} \cdot \lambda$, for all $k$. In particular that  $f_{i_k}^{p_k}$ must perform the reflection $s_{i_k}$ on the weight $s_{i_{k-1}} \ldots s_{i_1} \cdot \lambda .$ Therefore, 
$$p_k = \langle \alpha_{i_k}^\vee  , s_{i_{k-1}} \ldots s_{i_1} \cdot \lambda \rangle
= \langle w_{k-1}^{\bf i} \alpha_{i_k}^\vee, \lambda \rangle.$$
\end{proof}

\begin{Lemma} Let $b_\lambda \otimes c$ be a highest weight element of $B_\lambda \otimes B_\mu$. 
Let $b \otimes b_\mu^{\text{low}}$ be the lowest weight element of the component containing 
$b_\lambda \otimes c$. Let $ (p_1, \ldots p_m)$ be the downward Kashiwara data for $c$ with respect to ${\bf i}$, and $(q_1, \ldots q_m)$ the upward Kashiwara data for $b$ with respect to ${\bf i}^\rev := (i_m, \dots, i_1) $. Then, for all k, 
$p_k + q_{m-k+1} = \langle w_{k-1}^\wi \cdot \alpha_{i_k}^\vee ,  \nu\rangle,$ where $\nu= \wt (b_\lambda \otimes b)$.
\label{updowndata}
\end{Lemma}

\begin{proof}
Let $r_k= \langle w_{k-1}^\wi \cdot \alpha_{i_k}^\vee , \nu\rangle.$ By part (i) of Lemma \ref{KDL}, $(r_1, \ldots r_m)$ is the downward Kashiwara data for $b_\lambda \otimes c$. Define $b_k \in B_\lambda$ and $c_k \in B_\mu$ by 
$b_k \otimes c_k = f_{i_k}^{r_k} \ldots f_{i_1}^{r_1} ( b_\lambda \otimes c )$. Part (ii) of Lemma \ref{KDL}, along with the definition of Kashiwara data, shows that, for each $1 \leq k \leq m$, 
$$e_{i_k} (b_{k-1} \otimes c_{k-1})= 0 \hspace{0.2in} \mbox{and} \hspace{0.2in} f_{i_k} (b_k \otimes c_k) = 0.$$
In particular, the tensor product rule for crystals implies
$$e_{i_k} b_{k-1} = 0 \hspace{0.3in} \mbox{and} \hspace{0.2in} 
f_{i_k} c_k  = 0.$$
Define $p_k$ to be the number of times $f_{i_k}$ acts on $c_{k-1}$ to go from $b_{k-1} \otimes c_{k-1}$ to
$b_k \otimes c_k$, and $q_{m-k+1}$ to be the number of times $f_{i_k}$ acts on $b_{k-1}$. Since $f_{i_k} c_k=0$, we see that $\varphi_{i_k} (c_{k-1}) = p_k$. Hence, by definition $(p_1, \ldots p_m)$ is the downward Kashiwara data for $c$ with respect to $ \wi $. Similarly,
$e_{i_k} b_{k-1} = 0$, so $\varepsilon_{i_k} (c_k) = q_{m-k+1}$. By Lemma \ref{th:bottom}, $b_m= b$, so this implies that $(q_1, \ldots q_m)$ is the upward Kashiwara
data for $b$ with respect to $ \wi^\rev $. 
Since $p_k + q_{m-k+1}$ is the number of times that $ f_{i_k} $ acts on $ b_{k-1} \otimes c_{k-1} $ to reach $b_k \otimes c_k$, we see that $ p_k + q_{m-k+1} = r_k $.
\end{proof}

Let $b_\lambda \otimes c$ be a highest weight element in $B_\lambda \otimes B_\mu$. As discussed in section \ref{maintheorem}, $*c$ can be considered as an element of $B_\lambda$.

\begin{Lemma} Define $\nu = wt (b_\lambda \otimes c)$. Let $(p_1, \ldots p_m)$ be the downward Kashiwara data for $c \in B_\mu$ with respect to ${\bf i}.$ Let  $(q_1, \ldots q_m)$ be the downward Kashiwara data for $*c \in B_\lambda$ with respect to the decomposition 
$\theta ({\bf i}^\rev) := (\theta(i_m), \dots, \theta(i_1))$ of $w_0$. Then, for all $k$, 
$p_k +q_{m-k+1} = \langle w_{k-1}^\wi \cdot \alpha_{i_k}^\vee, \nu \rangle.$
\label{*data}
\end{Lemma}

\begin{proof}
The proof will depend on results from \cite{Kam} on the MV polytope model for crystals.  In particular, within this model it is easy to express Kashiwara data and the Kashiwara involution.

Let $ P = P(M_\bigdot) $ be the MV polytope of weight $ (\nu - \lambda, \mu) $ corresponding to $ c $.  Then by Theorem 6.6 of \cite{Kam}, 
\begin{equation*} \label{eq:pM} 
p_k = M_{w_{k-1}^\wi \cdot \Lambda_{i_k}} - M_{w_k^\wi \cdot \Lambda_{i_k}}.
\end{equation*}

Now, consider $ P $ as a stable MV polytope (recall that this means that we only consider it up to translation).  Then by Theorem 6.2 of \cite{Kam}, we see that $ *(P) = -P $.

The element $ \iota_\mu^{-1} *\iota_\lambda(c) \in B_\lambda $ corresponds to the MV polytope $ \nu - P $ and hence has BZ datum $ N_\bigdot $,  where $M_\bigdot$ and $N_\bigdot$ are related by
\begin{equation*} \label{eq:N}
M_\gamma = \langle \gamma, \nu \rangle + N_{-\gamma}.
\end{equation*}

Let $\wi' = \theta(\wi^{\rev}) $.  
Then,
\begin{equation*}
 -w_{k-1}^{\wi} \cdot \Lambda_{i_k} = w_{m-k +1}^{\wi'} \cdot \Lambda_{i'_{m-k+1}} \text{ and } 
 -w_k^{\wi} \cdot \Lambda_{i_k} = w_{m-k}^{\wi'} \cdot \Lambda_{i'_{m-k+1}}.
 \end{equation*}

Combining the last 3 equations, we see that
\begin{equation} \label{eqpk}
\begin{aligned}
p_k &= \langle w_{k-1}^{\wi} \cdot \Lambda_{i_k}, \nu \rangle + N_{-w_{k-1}^{\wi} \cdot \Lambda_{i_k}} - \langle w_k^{\wi} \cdot \Lambda_{i_k}, \nu \rangle - N_{-w_k^{\wi} \cdot \Lambda_{i_k}} \\
&= N_{w_{m-k+1}^{\wi'} \cdot \Lambda_{i'_{m-k+1}}} - N_{w_{m-k}^{\wi'} \cdot \Lambda_{i'_{m-k+1}}} + \langle w_{k-1}^{\wi} \cdot \Lambda_{i_k} - w_k^{\wi} \cdot \Lambda_{i_k}, \nu \rangle.
\end{aligned}
\end{equation}

Applying Theorem 6.6 of \cite{Kam} again,
\begin{equation} \label{eqqk}
q_k = N_{w_{k-1}^{\wi'} \cdot \Lambda_{i'_k}} - N_{w_k^{\wi'} \cdot \Lambda_{i'_k}}.
\end{equation}

We now add equation (\ref{eqpk}) and (\ref{eqqk}), substituting $m-k+1$ for $k$ in the second equation, to get
\begin{equation*} \label{eq:p+q} p_k+q_{m-k+1}= \langle w_{k-1}^{\wi} \cdot \Lambda_{i_k} - w_k^{\wi} \cdot \Lambda_{i_k}, \nu \rangle = \langle w_{k-1}^{\wi} \cdot (\Lambda_{i_k} - s_{i_k} \cdot \Lambda_{i_k}), \nu \rangle = \langle w_{k-1}^{\wi} \cdot \alpha_{i_k}^\vee, \nu \rangle.
\end{equation*}
\end{proof}

\begin{proof}[{\bf Proof of Theorem \ref{maintheorem}}]
We know that $\sigma_{B_\lambda, B_\mu} (b_\lambda \otimes c) = b_\mu \otimes b$ for some $b \in B_\lambda$. By the definition of $\sigma_{B_\lambda, B_\mu}$ (see  Section \ref{commutor_section}),  we see that 
\begin{equation*}
\xi (b_\lambda \otimes c)=  (\xi \circ \xi) (\mathrm{Flip}(\sigma(b_\lambda \otimes c))) =  \xi(b) \otimes b_\mu^{low}.
\end{equation*}
In particular, $\xi(b)  \otimes b^{low}_\mu$ is the lowest weight element of the component of $B_\lambda \otimes B_\mu$ containing $b_\lambda \otimes c$. 

Fix a reduced expression ${\bf i}= (i_1, \ldots, i_m)$ for $w_0$. Let $(p_1, \ldots, p_m)$ be the downward Kashiwara data for $c$ with respect to ${\bf i}$, and let $(q_1, \ldots q_m)$ be the downward Kashiwara data for $b$ with respect to $\theta({\bf i}^\rev):= (\theta(i_m), \cdots, \theta(i_1))$. 
Notice that $(q_1, \ldots q_m)$ is also the upward Kashiwara data for $\xi (b)$ with respect to ${\bf i}^{\rev}:= (i_m, \ldots i_1)$, since $\xi$ interchanges the action of $f_i$ and $e_{\theta(i)}$. Hence by Lemma \ref{updowndata}, 
\begin{equation} \label{sequ}  p_k+ q_{m-k+1} = \langle w_{k-1}^\wi \cdot \alpha_{i_k}^\vee, \nu \rangle \end{equation}
for all $k$, where $\nu$ is the weight of $b_\lambda \otimes c$.

As discussed in Section \ref{mainsection}, $*c \in \iota(B_\mu)$, and so can be considered as an element of $B_\mu$. Let $(q_1', \ldots q_m')$ by the downward Kashiwara data for $*c \in B_\mu$ with respect to $\theta ( {\bf i}^{\rev} )$. 
By Lemma \ref{*data} we have 
\begin{equation} \label{*equ} p_k + q_{m-k+1}' = \langle w_{k-1}^\wi \cdot \alpha_{i_k}^\vee, \nu \rangle . \end{equation}
Comparing equations (\ref{sequ}) and (\ref{*equ}), $q_k = q_k'$ for each $1 \leq k \leq m$. That is, the downward Kashiwara data for $b$ and $*c$ with respect to $\theta({\bf i}^{\rev})$ are identical.
Hence by Lemma \ref{th:bottom}, $ b = *c $. \end{proof}

\section{Questions} \label{questions}

The involution $ * $ gives $ B_\infty $ an additional crystal structure, defined by $ f_i^* \cdot b := * \circ f_i \circ * (b) $.  Let $ B_\infty^i $ denote the crystal with vertex set $\mathbb{Z}_{\geq 0}$, where $ e_j, f_j $ act trivially for $ j \ne i $ and $ e_i, f_i $ act as they do on the usual $B_\infty $ for $ \mathfrak{sl}_2 $. Kashiwara \cite[Theorem 2.2.1]{Kas1} showed that the map
\begin{align*}
B_\infty &\rightarrow B_\infty \otimes B_\infty^i \\
b &\mapsto (e_i^*)^{\varepsilon_i(*b)}(b) \otimes \varepsilon_i(*b)
\end{align*}
is a morphism of crystals with respect to the usual crystal structures on each side.  We can think of this fact as an additional property of $ * $.

On the other hand, the commutor $ \sigma $ also has an additional property, which is called the cactus relation.  This relation states that if $ A, B, C $ are crystals, then $$ \sigma_{A, C \otimes B} \circ (1 \otimes \sigma_{B,C}) = \sigma_{B \otimes A, C} \circ (\sigma_{A, B} \otimes 1). $$ (See \cite[Theorem 3]{HK}).

\begin{Question}
Is there a relation between this additional property of Kashiwara's involution $ * $ and the cactus relation for the commutor $ \sigma$?
\end{Question}

Another direction is to consider the generalization beyond finite dimensional reductive Lie algebras.  We can define a crystal commutor for any the crystals of highest weight representations of any symmetrizable Kac-Moody algebra by $\sigma (b_\lambda \otimes c) = b_\mu \otimes *c$ whenever $ b_\lambda \otimes c $ is a highest weight element. This will be well defined by the analysis given in section \ref{mainsection}.
\begin{Question}
Does this commutor satisfy the cactus relation? \footnote{This question has recently been answered in the affirmative by Savage \cite[Theorem 6.4]{S}}
\end{Question}

\end{document}